\documentclass{article}
\usepackage[utf8]{inputenc}

\usepackage{float} 
\usepackage{amsmath}
\usepackage{amssymb}
\usepackage{array}
\usepackage{bbm}
\usepackage{makeidx}
\usepackage{hyperref}
\usepackage{cleveref}
\usepackage{caption}
\usepackage{subcaption}
\usepackage{amsthm}
\usepackage{color}
\usepackage{mathtools}
\usepackage{graphicx} 
\usepackage[LGR,T1]{fontenc}
\usepackage{siunitx}
\usepackage{mathrsfs}
\usepackage{scalerel}
\usepackage{forest}
\usepackage{algorithm,algpseudocode}
\usepackage{lipsum}
\usepackage{booktabs}
\usepackage{siunitx}
\usepackage{xcolor}
\usepackage[title]{appendix}
\hypersetup{
    colorlinks,
    linkcolor={blue!80!black},
    citecolor={green!50!black},
}
\usepackage{accents}
\usepackage{apptools}
\usepackage{tikz}
\usetikzlibrary{arrows, positioning, automata}
\usepackage{comment}
\usepackage[shortlabels]{enumitem}

\usepackage[mathscr]{euscript}
\DeclareSymbolFont{rsfs}{U}{rsfs}{m}{n}
\DeclareSymbolFontAlphabet{\mathscrsfs}{rsfs}

\def\B{{\rm Band}}
\def\TAP{{\rm TAP}}
\def\Chull{\overline\Sigma_N}

\def\cuU{\mathscrsfs{U}}
\def\cuL{\mathscrsfs{L}}
\def\de{{\rm d}}
\def\Par{{\sf P}}
\def\sTV{\mbox{\tiny \rm TV}}

\def\<{\langle}
\def\>{\rangle}

\def\sEA{\mbox{\tiny \rm EA}}
\def\sTV{\mbox{\tiny \rm TV}}

\def\salg{\mbox{\tiny \rm alg}}
\def\sEA{\mbox{\tiny \rm EA}}

\def\prob{{\mathbb P}}

\AtAppendix{\counterwithin{lemma}{section}}

\usepackage{mathtools}

\newtheoremstyle{myremark} 
    {\topsep}                    
    {\topsep}                    
    {\rm}                        
    {}                           
    {\bf}                        
    {.}                          
    {.5em}                       
    {}  

\newtheorem{claim}{Claim}[section]

\newtheorem{theorem}{Theorem}

\newtheorem{corollary}[claim]{Corollary}
\newtheorem{definition}[claim]{Definition}

\theoremstyle{myremark}
\newtheorem{remark}{Remark}[section]

\allowdisplaybreaks
\usepackage[top=1in,bottom=1in,left=1in,right=1in]{geometry}
\usepackage[utf8]{inputenc}

\def\ed{\stackrel{{\mathrm d}}{=}}

\def\bA{{\boldsymbol A}}

\def\cuP{\mathscrsfs{P}}

\def\cA{{\mathcal A}}

\def\cS{{\mathcal S}}

\def\cT{\mathcal{T}}
\def\cV{\mathcal{V}}

\def\cS{\mathcal{S}}

\def\GOE{\mbox{GOE}}

\def\Vol{{\rm Vol}}

\def\hbm{\hat{\boldsymbol m}}

\def\bsigma{{\boldsymbol \sigma}}
\def\S{{\mathbb S}}

\newcommand{\Crt}{\mathop{\mathrm{Crt}}\nolimits}
\newcommand{\indi}{\ensuremath{\boldsymbol 1}}

\def\bA{\mathbf{A}}

\def\bG{\mathbf{G}}

\def\bM{\mathbf{M}}

\def\bP{\mathbf{P}}

\def\bW{\mathbf{W}}

\def\bm{\boldsymbol{m}}

\def\bv{\boldsymbol{v}}

\def\bx{\boldsymbol{x}}

\def\bz{\boldsymbol{z}}

\def\bA{\boldsymbol{A}}

\def\bG{\boldsymbol{G}}

\def\bM{\boldsymbol{M}}

\def\bP{\boldsymbol{P}}

\def\bW{\boldsymbol{W}}

\def\normal{{\mathsf{N}}}

\def\bzero{\boldsymbol{0}}
\def\bfone{\boldsymbol{1}}

\def\bm{\boldsymbol{m}}

\def\oSigma{\overline{\Sigma}}
\def\conv{{\rm conv}}
\def\supp{{\rm supp}}

\def\EG{{\sf E}}
\def\THR{{\sf THR}}
\def\ALG{{\sf ALG}}
\def\OPT{{\sf OPT}}
\def\SOC{{\sf SOC}}
\def\Lev{{\sf L}}
\def\Hess{\boldsymbol{\mathcal H}}

\def\reals{{\mathbb R}}
\def\naturals{{\mathbb N}}

\def\sT{{\sf T}}

\def\z{{\boldsymbol z}}
\def\u{{\boldsymbol u}}
\def\Ball{{\sf B}}
\def\Round{{\sf Round}}

\def\id{{\boldsymbol I}}

\def\rP{{\rm P}}

\renewcommand{\P}{\mathbb{P}}
\newcommand{\E}{\mathbb{E}}
\newcommand{\R}{\mathbb{R}}

\newcommand{\eps}{\varepsilon}

\newcommand{\Unif}{\operatorname{Unif}}

\newcommand{\RN}[1]{%
  \textup{\uppercase\expandafter{\romannumeral#1}}%
}

\newcommand{\RNum}[1]{\uppercase\expandafter{\romannumeral #1\relax}}



\makeatletter
\newcommand*{\rom}[1]{\expandafter\@slowromancap\romannumeral #1@}
\makeatother

\title{Optimization of random high-dimensional functions:\\
Structure and algorithms}

\author{Antonio Auffinger\thanks{Department of Mathematics, Norhtwestern University}
	\and 
	Andrea Montanari
	\thanks{Department of Electrical Engineering and Department of Statistics, Stanford University} 
	\and 
	Eliran Subag\thanks{Department of Mathematics, Weizmann Institute of Science}
}
\date{}
\begin{document}
\maketitle

\begin{abstract}
Replica symmetry breaking postulates that near optima of 
spin glass Hamiltonians have an ultrametric structure. Namely, near optima can be 
associated to leaves of a tree, and the Euclidean distance between them corresponds
to the distance along this tree. 
We survey recent progress towards a rigorous proof of this picture
in the context of  mixed $p$-spin spin glass models.
We focus in particular on the following topics: $(i)$~The structure of critical points
of the Hamiltonian; $(ii)$~The realization of the ultrametric tree as near optima of
a suitable TAP free energy; $(iii)$~The construction of efficient optimization
algorithm that exploits this picture.
\end{abstract}

\section{Introduction}

Mean field spin glasses are high-dimensional random functions $H_N:\Sigma_N\to\reals$,
$\Sigma_N\subseteq \reals^N$ with special symmetry properties.  
In the most studied cases, for any $k$ points (spin configurations) 
$\bsigma_1,\dots,\bsigma_k\in\Sigma_N$,
the joint distribution of $H_N(\bsigma_1)$, \dots $H_N(\bsigma_k)$
only depends on the configurations through their joint empirical distribution
$N^{-1}\sum_{i=1}^N\delta_{\sigma_{1,i},\dots,\sigma_{k,i}}$. 

These models were originally motivated by the study of disordered magnetic materials.
However, it soon become clear that a large number of random optimization problems of interest in
computer science and statistics fit this framework \cite{mezard1987spin}. In this article
we survey recent rigorous work aimed at describing the structure of near optima of
spin glasses, with a focus onto questions that are relevant for optimization.

We will consider the most classical mean-field spin glass model, namely the mixed 
$p$-spin spin glass, with either $\Sigma_N=\{+1,-1\}^N$ (Ising model) or 
$\Sigma_N=\S^{N-1}(\sqrt{N})$ (spherical model). The Hamiltonian is defined by
\begin{align}
H_N(\bsigma) & := \sum_{k\ge 2}\frac{c_k}{N^{(k-1)/2}}\<\bG^{(k)},\bsigma^{\otimes k}\>
\label{eq:FirstHamiltonian}\\
& = \sum_{k\ge 2}\frac{c_k}{N^{(k-1)/2}}\sum_{i_1,\dots,i_k=1}^N 
G_{i_1\dots i_k}^{(k)}\sigma_{i_1}\cdots\sigma_{i_k}\, ,
\end{align}
where $\bG^{(k)}=(G^{(k)}_{i_1\dots i_k})_{i_1,\dots,i_k\le n}$ is a tensor
with i.i.d. Gaussian components $G^{(k)}_{i_1\dots i_k}\sim\normal(0,1)$.

It is immediate to see that $\{H_N(\bsigma)\}_{\bsigma\in \Sigma_N}$ is a centered Gaussian 
process on $\Sigma_N$,
with covariance
\begin{align}
\E\big[H_N(\bsigma)H_N(\bsigma')\big] = N \xi\big(\langle\bsigma,\bsigma'\rangle/N\big)\, ,
\end{align}
where $\xi(x) := \sum_{k\ge 2}c_k^2x^k$. A term linear in $\bsigma$ 
(a `magnetic field' in physics language) could be added but we omit it here because this simplifies 
some of the statements and discussions below. We will assume throughout that 
the non-random coefficients $c_k$ decay fast enough that  $\sum_{k\ge 2}c_k^2k^{C}<\infty$ for any $C>0$.

Understanding the structure of optima and near optima of spin glass Hamiltonians
is a central concern of the theory. The most basic quantity is of course the optimum 
value
\begin{align}
\OPT_N:= \frac{1}{N}\E\max_{\bsigma\in \Sigma_N}H_N(\bsigma)\, ,\;\;\;\; \OPT :=
\lim_{N\to\infty} \OPT_N\, . 
\end{align}
(The standard physics convention is to study the minimum of $-H_N(\bsigma)$,
but of course the two formulations are equivalent.)
We occasionally use $\OPT_N(H_N) := \max_{\bsigma\in \Sigma_N}H_N(\bsigma)/N$ 
for the maximum of the random Hamiltonian $H_N$.
The $N\to\infty$ limit above was proven to exist and to be given by a `zero-temperature' 
version of Parisi's formula, which is a variational principle over the following space of functions:
\begin{align}\label{def:space}
\cuU :=\Big\{\gamma:[0,1) \to\reals_{\ge 0}: \;\; \gamma\mbox{ non-decreasing }, \int_{0}^1\gamma(t)\,\de t<\infty\Big\} \, .
\end{align}
For $\gamma\in \cuU$, let $\Phi_{\gamma}:[0,1]\times \reals\to \reals$ be the solution of the following PDE, known as \emph{the Parisi PDE}
\begin{align}\label{eq:PDEFirst}
\partial_t \Phi_{\gamma}(t,x)+\frac{1}{2}\xi''(t) \Big(\partial_x^2\Phi_{\gamma}(t,x)+\gamma(t) (\partial_x\Phi_{\gamma}(t,x))^2\Big) = 0\, ,
\end{align}
with terminal condition at $t=1$:
\begin{align}
\Phi_{\gamma}(1,x) = \begin{cases}
x^2/2 & \mbox{if $\Sigma_N=\S^{N-1}(\sqrt{N})$,}\\
|x|& \mbox{if $\Sigma_N=\{+1,-1\}^N$.}\\
\end{cases}
\end{align}
The Parisi functional $\Par:\cuU\to \reals$ is then defined by
\begin{align}\label{eq:ParisiF}
\Par(\gamma) \equiv \Phi_{\gamma}(0,0)-\frac{1}{2}\int_{0}^1 t\xi''(t)\gamma(t)\, \de t\, .
\end{align}
The zero-temperature variational principle is stated below: its proof builds on
earlier results at positive temperature proved by Talagrand \cite{talagrand2006parisi} and
Panchenko \cite{panchenko2013parisi,panchenko2013sherrington}. 
\begin{theorem}[\cite{ArnabChen15, JT,auffinger2017parisi}]\label{thm:ZeroTParisi}
The limit $\OPT :=
\lim_{N\to\infty} \OPT_N$ exists and is given by
\begin{align}
\OPT = \inf_{\gamma\in\cuU} \Par(\gamma)\, .\label{eq:FirstParisi}
\end{align}
\end{theorem} 
We also recall that classical Gaussian concentration results imply that
$\max_{\bsigma\in\Sigma_N}H_N(\bsigma)/N$ is exponentially concentrated around $\OPT_N$
and therefore the asymptotic value $\OPT$ is representative of the typical optimal value
for  large $N$.

Of course, characterizing the typical value of the optimum is only a first step 
towards understanding the landscape structure.
A  more detailed picture is obtained by considering the random 
superlevel sets
\begin{equation}\label{sublevel.sets}
\Lev_{N}(\eta):= \big\{ \bsigma \in \Sigma_N : H_N(\bsigma) \ge N \eta \big\},
\end{equation}
where $\eta$ is a fixed real number. 
It is clear that for $\eta_{1} \le \eta_{2},$ $\Lev_{N}(\eta_{2}) \subseteq \Lev_{N}(\eta_{1})$. 
Of course, the maximum of $H_{N}$  can be written as 
\begin{align}
\OPT_N(H_N) := \frac{1}{N}\max_{\bsigma\in \Sigma_N}H_N(\bsigma)
  = \sup \big\{ \eta \in \mathbb R : \Lev_{N}(\eta)\neq \emptyset  \big\}\, .
\end{align}

An interesting approach to study the geometry of the random set 
$\Lev_N(\eta)$ is to consider the uniform measure over this set 
$U_{N,\eta}(\de\bsigma) \propto \bfone_{\Lev_N(\eta)}(\bsigma)\de\bsigma$. This turns out
to be closely approximated by the following Gibbs measure
(for a suitable choice of $\beta = \beta_*(\eta)$)
\begin{align}
G_{N,\beta}(\de \bsigma) = \frac{1}{Z_{N,\beta}} e^{\beta H(\bsigma)}\, \de\bsigma\, .
\label{eq:Gibbs}
\end{align}
Here $\de\bsigma$ is the uniform measure over the sphere
$\Sigma_N=\S^{N-1}(\sqrt{N})$ or the 
counting measure over the hypercube
$\Sigma_N=\{+1,-1\}^N$ (in the latter case we keep writing integrals instead of sums.)

The free energy density is defined as the exponential growth rate
of partition function $Z_{N,\beta}$:
\begin{equation}
	F_{N,\beta} := \frac1N\log Z_{N,\beta}=\frac1N\log\int_{\Sigma_N}e^{\beta H_N(\bsigma)}\, 
	\de\bsigma\, . \label{eq:FirstFree}
\end{equation}
Roughly speaking (and if $\beta=\beta_*(\eta)$ is chosen as mentioned above),
we have $F_{N,\beta}\approx \beta\eta+N^{-1}\log {\rm Vol}(\Lev_N(\eta))$.

The rest of this article is organized as follows:
\begin{itemize}
\item We begin by discussing some fundamental geometric and topological
properties of the optimization landscape in 
Section \ref{sec:Critical}. In particular, the number of local
maxima of $H_N(\bsigma)$ in the superlevel set $\Lev_N(\eta)$
is often exponentially large in the dimension $N$. Further, nearly orthogonal 
vectors can be found in $\Lev_N(\eta)$, for $\eta$ close to $\OPT$.
\item While the previous section indicates that the spin glass landscape is not only
non-convex  but in fact extremely rough,
Section \ref{sec:Free} outlines the hidden tree structure behind this seemingly 
unstructured landscape. We switch from the level sets $\Lev_{N}(\eta)$
to the Gibbs measure $G_{N,\beta}(\de \bsigma)$. The space $\Sigma_N$ can be partitioned
hierarchically into sets (`ancestor states') of decreasing normalized radius (from $1$ to 
$(1-q_{\sEA})^{1/2}$), and to each of these sets $\alpha$ can be associated its 
barycenter\footnote{Here and below $\conv(S)$ denotes the convex hull of a set $S\in\reals^d$.}
$\bm^{(\alpha)}\in\oSigma_N:=\conv(\Sigma_N)$. The vectors $\bm^{(\alpha)}$ are organized 
according to an ultrametric tree.
\item Can we use the ultrametric tree to explore the structure of optima or near optima
of $H_N(\bsigma)$? Section \ref{sec:algorithms} shows that this is indeed the case by constructing
two efficient algorithms that output configurations $\bsigma^{\salg}\in\Sigma_N$
by constructing random paths in $\oSigma_N$. The asymptotic value achieved by these algorithms 
$\ALG := \lim_{N\to\infty}H_N(\bsigma^{\salg})/N$ turns out to be given by a modified
Parisi formula. In particular, when the optimizer of Eq.~\eqref{eq:FirstParisi} is
strictly monotone (`no overlap gap'), we have $\ALG = \OPT$.
\end{itemize}

\section{Critical points and the landscape structure}
\label{sec:Critical}

In this section, we focus on two fundamental questions:
\begin{enumerate}
\item[$(a)$]  Computing the number of critical points 
(local maxima, saddles) of $H_{N}$ at different values of the energy.  
In particular, we describe the exponential growth rate of this number, known in physics as 
the complexity function.
\item[$(b)$] Studying the geometry of the  superlevel set $\Lev_N(\eta)$ when
$\eta$ is close to the maximum $\OPT$. In particular, we will discuss the connection
between Parisi's replica symmetry breaking formula and the
range of values of $\<\bsigma^1,\bsigma^2\>/N$ when $\bsigma^1,\bsigma^2\in\Lev_N(\eta)$.
\end{enumerate}

Throughout this section we will refer to  $\xi$ as a mixture and use the shorthands: 
\begin{equation*}
\xi:= \xi(1)>0\, ,\;\;\;\;  \xi':= \xi'(1) > 0\, ,\;\;\;  \xi'':= \xi''(1) >0. 
\end{equation*}
We say that the model is pure (of degree $p$) if $\xi(x) = c_p^2x^p$, which amounts to say that the 
Hamiltonian is an homogeneous polynomial of degree $p$.
Note that $\xi''\geq \xi'$ with equality only in the pure case with $p=2$. 

\subsection{Complexity of critical points}

Consider the case where $\Sigma_N=\S^{N-1}(\sqrt{N})$ and the subsets \eqref{sublevel.sets}
 are smooth submanifolds of $\S^{N-1}(\sqrt{N})$.  We introduce the complexity of spherical 
 spin glasses as follows. For any  $\eta \in \R$ and integer $0\le k < N$, we consider the
  (random) number $\Crt_{N,k}(\eta)$ of critical points of the
Hamiltonian $H_N$ in the set $\Lev_{N}(\eta)$  with index equal to
$N-1-k$,
\begin{equation}
  \label{defWk}
  \Crt_{N,k}(\eta) = 
  \sum_{\boldsymbol\sigma: \nabla H_N(\boldsymbol\sigma) = 0 } 
  \indi\big\{ \bsigma \in \Lev_{N}(\eta)\big\}\; \indi\big\{
    i(\nabla^2 H_N(\boldsymbol\sigma)) = N-1-k\big\}.
\end{equation}
Here $i(\nabla^2 H_N(\boldsymbol\sigma))$ is the index of 
$\nabla^2 H_{N}$ at $\boldsymbol\sigma$, that is, the number of negative eigenvalues 
of the Hessian $\nabla^2 H_{N}(\bsigma)$.  As an illustration, $\Crt_{N,0}(\eta)$ counts 
the number of local maxima with energy value larger than $N\eta$, while $\Crt_{N,1}(\eta)$ 
counts the number of saddles in $\Lev_{N}(\eta)$ with ``one positive direction,''
and so on. 
 Here $\nabla$, $\nabla^2$ are the Riemannian gradient and Riemannian Hessian which act on
the tangent space of $\S^{N-1}(\sqrt N)$.

We will also consider the (random) total number 
$\Crt_{N}(\eta)$ of critical values of the Hamiltonian $H_N$ in the set $\Lev_{N}(\eta)$
(whatever their index)
\begin{equation}
  \label{e:defW}
  \Crt_{N}(\eta) 
  = \sum_{\boldsymbol\sigma: \nabla H_N(\boldsymbol\sigma) = 0 } 
  \indi\big\{ \bsigma \in L_{N}(\eta)\big\} .
\end{equation}
Note that $\Crt_{N,k}(\eta) \leq \Crt_{N}(\eta)$ for all $k\geq 0$.

For general mixtures, the asymptotics of the mean of $\Crt_{N,k}(\eta)$ was first derived in
\cite{AB}. 
\begin{theorem}[Theorem 1 in \cite{AB}]\label{thmAB} For any $\eta \in \mathbb R$, and any fixed $k \geq 0$, there exist $\Sigma_{k}(\eta), \Sigma(\eta) \in \mathbb R$ such that
\begin{equation*}
    \lim_{N\to\infty}
    \frac{1}{N} \log \E \Crt_{N,k}(\eta) = \Sigma_{k}(\eta),
\end{equation*}
and
   \begin{equation*}
    \lim_{N\to\infty}
    \frac{1}{N} \log \E \Crt_{N}(\eta) = \Sigma(\eta).
\end{equation*}
\end{theorem}
The exact expression for $\Sigma(\eta),  \Sigma_{k}(\eta)$ can be found in \cite[Equations (2.13)]{AB}. These functions
 depend on the model $\xi$ (we avoid explicitly writing this dependence in our notation) and are
  called the (averaged or annealed) complexity of critical points and (averaged or annealed) 
  complexity of critical points of index $N-1-k$.

  The complexity functions $\Sigma_k$, $\Sigma$ have a number of interesting properties, 
  which we next summarize.
Set 
\begin{equation}
 \THR' = \frac{2 \xi' \sqrt{\xi''}}{(\xi'+\xi'')\xi}, \quad\quad \THR = \frac{ (\xi'' - \xi')\xi + \xi'^2 }{\xi'\sqrt{\xi\xi''}}.
\end{equation}
Then we have the following:
\begin{enumerate}
\item[$(i)$] $\THR'\leq \THR$ with equality if and only the model is pure.
\item[$(ii)$] The functions $\eta \to \Sigma_{k}(\eta)$  and  $\eta \to \Sigma(\eta)$ are continuous and non-increasing.
(Note that monotonicity follows by definition.)
\item[$(iii)$] For each $k$, $\Sigma_{k}(\eta)$ is strictly decreasing on $(\THR',\infty)$, with 
\[
\lim_{\eta \to \infty}\Sigma_{k}(\eta) = -\infty.
\] 
Further, on $(-\infty,\THR')$, $\Sigma_k$ is constant (and taking its maximum value) with
\begin{equation*}\label{totalmaxima}
 \eta\in (-\infty,\THR') \;\;\Rightarrow\;\; \Sigma_k(\eta)=\frac{1}{2} \log \frac{\xi''}{\xi'}-\frac{\xi''-\xi'}{\xi''+ \xi'} > 0.
\end{equation*}
\item[$(iv)$] For any $k,k'\in \mathbb N$ with $k < k'$, $\Sigma_{k}(\eta) > \Sigma_{k'}(\eta) $ for all 
$\eta \in(\THR,\infty)$.  
\item[(iv)] For $\eta\in (\THR',\infty)$, $\Sigma_{0}(\eta) = \Sigma(\eta)$.

\end{enumerate}

What do these properties and Theorem \ref{thmAB} above tell us about the landscape? First, define 
$\EG_{k} \in \mathbb R$ to be the unique energy value such that 
\begin{align}
\Sigma_{k} (\EG_{k}) = 0\, . 
\end{align}
These quantities are well-defined due to properties $(ii)$ and $(iii)$.  A simple application of Markov's inequality yields
\[ \mathbb P \left(\OPT_N(H_N) \ge\eta \right) = \mathbb P (\Lev_{N}(\eta) \neq \emptyset) = 
 \mathbb P (\Crt_{0}(\eta) \geq 1) \leq \mathbb E \Crt_{0}(\eta) \,.
\]
Combined with Theorem \ref{thmAB} this  leads to the bound 
\begin{equation}\label{osb}
\OPT \leq \EG_{0}.
\end{equation}
Note that, in general, we cannot argue in the opposite direction. Even if $\Sigma_0(\eta)>0$,
we cannot conclude that with high probability there exist local maxima in $\Lev_N(\eta)$. 
However, as discussed below, in some cases this conclusion holds.

In fact, Markov's inequality also implies that, for any $\varepsilon >0$, the event that there is 
a critical value of the Hamiltonian $H_{N}$ above 
  level $N(\EG_{k}+\varepsilon)$ and with index $N-1-k$ (or larger) has an exponential small probability.
   Loosely speaking, as $N$ diverges, there are no saddles with $k$ positive directions 
   above $\EG_{k}$. This is a good moment to remind the reader that the scaling of
    $\Lev_N(\eta)$ contains a factor $N$: namely, $\Lev_N(\eta)$ is the set of
    spin configurations $\bsigma$ satisfying $H_N(\bsigma)\ge N\eta$. 
    If $\eta>\EG_1$, $\Lev_N(\eta)$ contains no saddles, with high probability (by argument given above). 
    If $\Lev_N(\eta)$ is non-empty and contains multiple local maxima,
    it must be the case that  $\Lev_N(\eta)$ is formed by disconnected components. 
    Each component contains one local maximum and they can only be connected if we decrease the 
    energy by an amount $\Theta(N)$. In other words, in such setting, energy barriers are of order $N$. 
    As we will see shortly, this scenario indeed holds for the spherical pure $p$-spin.

First of all, 
the formulas for the complexity functions $\Sigma_{k}(\eta), \Sigma(\eta)$ have a simpler 
form for the spherical pure $p$-spin. If this is the case, we can assume, without loss of generality, $\xi(t)=t^{p}$. 
It was shown in \cite{ABC} that 
\begin{align}
\Sigma(\eta)=\begin{cases}
\frac{1}{2}\log(p-1) & \text{if }\eta\le 0,\\
\frac{1}{2}\log(p-1)-\frac{p-2}{4(p-1)}\eta^{2} & \text{if } 0\le \eta \leq\THR,\\
\frac{1}{2}\log(p-1)-\frac{p-2}{4(p-1)}\eta^{2}-J(\eta) & \text{if }\THR\le \eta,\\
\end{cases}
\end{align}
where $\THR=\THR(p)=2\sqrt{\frac{p-1}{p}}$ and for $\eta\ge\THR $,
\[
J(\eta)=\frac{\eta}{\THR^{2}}\sqrt{\eta^{2}-\THR^{2}}-\log\big(\eta+\sqrt{\eta^{2}-\THR^{2}}\big)+\log \THR.
\]
Further, in this case, $\THR=\THR'$.
A straightforward consequence of these formulas
(and of the strict monotonicity of $\Sigma_k(\eta)$ in $k$ for $\eta>\THR$ by property $(iv)$ above) 
is that in the spherical pure $p$-spin,  
\begin{equation}\label{eq:cascade}
\EG_{k} > \EG_{k+1} > \THR = \THR' = 2\sqrt{\frac{p-1}{p}} \quad \text{ for all } k \geq 0.
\end{equation}

Second, and more importantly, in the spherical pure $p$-spin the 
exponential growth rate of the typical number (as opposed to the average number)
of critical points can be rigorously computed.
This is referred to in physics as the \emph{quenched} complexity,
and turns out to coincide with the annealed complexity given above. This was established in \cite{subag2017complexity, SZ2021}.  

\begin{theorem}[Theorem 1 in \cite{subag2017complexity}] \label{Subag} For any $p\geq 3$ and $\eta \in(\THR,\,\EG_{0})$,
\begin{equation}
\lim_{N\to\infty}\frac{\mathbb{E}\big[\left(\Crt_N(\eta)\right)^{2}\big]}{\big(\mathbb{E}\left[ \Crt_N(\eta)\right] \big)^{2}}=1.\label{eq:moment_matching}
\end{equation}
Consequently, in $L^{2}$ and in probability,
\begin{equation}
  \label{eq-120521a}
\frac{\Crt_N(\eta)}{\E\Crt_N(\eta)} \overset{N\to\infty}{\longrightarrow} 1.
\end{equation}
\end{theorem}
The same result was also proved for $\eta \in(-\infty,\,\THR)$ in \cite{SZ2021} but only for $p\geq32$. It is expected to be true for any $p\geq3$. 
For models $\xi$ that are close to a pure $p$-spin, the logarithm of the ratio in \eqref{eq:moment_matching} was shown to be $o(N)$ for energies $\eta$ close to $\EG_{0}$, see \cite{benarous2020geometry}.  
Let us now explain in words some important consequences of Theorem \ref{Subag} above. First, one obtains 
that for any $\eta<\EG_0$, we have $\Crt_N(\eta) = \exp\{N\Sigma(\eta)+o(N)\}\ge 1$ with high probability.
This implies the inequality that  complements \eqref{osb}, namely
\[
\OPT\geq  \EG_{0},
\]
establishing, without the need use the Parisi Formula, the value of the ground state energy.

Further, the scenario described above of minima separated by barriers of order $N$
holds in this case by Theorem \ref{Subag} and Eq.~\eqref{eq:cascade}. Namely, for values of 
 $\eta \in (\EG_1, \EG_{1})$, $\Lev_N(\eta)$ can be written as 
\[
\Lev_N(\eta) = \bigcup_{\alpha = 1}^{M} C_{\alpha}
\]
where the sets $C_{\alpha}$ are connected, disjoint sets, each one homeomorphic to a point and 
$M =\exp(N \Sigma(\eta)+o(N))$. 

A further consequence of Theorem \ref{thmAB} --- always in the pure $p$-spin spherical model ---
 is the role played by the energy threshold $\THR$. This is the unique value of energy such that 
 the complexity of critical points of fixed index $k$ coincide, namely,
\[
 \Sigma_{k}(\THR) = \Sigma_{k'}(\THR) \text{ for all } k, k'\ge 0\, ,
\]
and such that for any $\varepsilon >0$ the probability that there exists a critical point of 
index $N-1-k$ in $H_N(\bsigma)\le N(\THR-\eps)$ decays exponentially for large
$N$, see \cite[Theorem 2.14]{ABC}. In other words, with high probability, for any fixed $k\geq 0$, a critical point of 
index $N-k-1$ can only exist at energy density values in the interval $[\THR,\EG_{k}]$. By symmetry, 
critical points of index $k$ can only exist at energy density values in the interval
$[-\EG_k,-\THR]$.

\begin{remark}[A short comment on methods]
The main tool to study the moments of the random variables $\Crt_{N}(\eta)$ and establish the 
results above is the Kac-Rice formula (see \cite{AT}, Chapter 11). This formula relates 
the $\ell$-th moment of $\Crt_{N}(\eta)$ to $\ell$-fold integrals of certain functions of 
$\ell$ determinants of $ \nabla^{2}H_{N} (\bsigma)$. 
The collection $\{\nabla^{2}H_{N} (\bsigma)\}_{\bsigma \in \Sigma_{N}}$ is a family of 
correlated random matrices whose entry distributions can be explicitly evaluated. It turns
 out that these Hessians can be analyzed as coupled Gaussian Orthogonal Ensembles plus a 
 random shift and tools from random matrix theory (in particular, large deviation principles)
  become available. We refer the reader to \cite{ABC}.
\end{remark}

\begin{remark}[Other directions]
Here we focused on the spherical mixed $p$-spin model. 
The complexity of the pure $p$-spin model was first computed using non-rigorous tools from 
statistical physics in \cite{crisanti1995TAP}.
The rigorous landscape-complexity program  --- counting critical 
points of high-dimensional random functions to understand their geometry --- initiated with a breakthrough insight by 
Y. Fyodorov in \cite{F1,F0}, who studied isotropic Gaussian fields followed by the work \cite{ABC}.
 More recently, this line of work led to the analysis of other disordered models including the bipartite spin glass 
 \cite{AChenBip, McKenna}, the spiked-tensor model \cite{BMMN19, PhysRevX.9.011003}, models with less invariance \cite{BPM21}, and Gaussian fields with isotropic increments 
 \cite{AZ20}.
\end{remark}

\subsection{The landscape near $\OPT$}

As we saw above, in the pure $p$-spin spherical model, the value of $\OPT$ can be 
computed from the annealed complexity of critical points $\Sigma(\eta)$ and it is equal to the 
unique point $\EG_{0}$ such that $\Sigma(\EG_0)=0$. In this section we also consider the 
Ising case where $\Sigma_{n} = \{ + 1,-1 \}^{N}$.

In the general case, the value of $\OPT$ is characterized by the zero-temperature
Parisi formula, cf. Theorem \ref{thm:ZeroTParisi}.
The minimizer $\gamma_{\star}$ of \eqref{eq:FirstParisi} is known to be unique, 
and is referred to as the Parisi measure at zero temperature. We write $\supp(\gamma_{\star})$ for its support, that is
\[
\supp (\gamma_{\star}) = \big\{ q \in [0,1] : \gamma_{\star} \text{ is not constant in any neighborhood of } q \big\}.
\]
The structure of  $\gamma_\star$ and in particular its support contains 
information about the structure of the superlevel set $\Lev_N(\eta)$
when $\eta$ is close to $\OPT$, as will be discussed below.

Physicists have put forward several predictions about the structure of $\gamma_\star$
for different models, on the basis of heuristic arguments and numerical solutions
of the variational principle. Unfortunately 
 verifying rigorously these predictions is very challenging, especially for the Ising
 model $\Sigma_{N}= \{+1,-1\}^{N}$, since in this case there is no closed form expression for $\gamma_\star$.

 In particular, in the case $\Sigma_{N}= \{+1, -1\}^{N}$, it is expected that there exists 
 $0\leq a<1$ such that  $[a,1] \subseteq  \mbox{supp } \gamma_{\star} $. This is known as 
 Full Replica Symmetry Breaking (FRSB) prediction at zero temperature.
 For the special case of the Sherrington-Kirkpatrick model, i.e. $\xi(t)=t^2$ and
 $\Sigma_N =\{+1,-1\}^N$, it is believed that $\supp(\gamma_\star)=[0,1]$.
  In this direction, the current state 
 of the art is given by the result of Auffinger-Chen-Zeng \cite{auffinger2017sk} that establishes
 that $\supp(\gamma_{\star} )$
  contains infinitely many points for any mixture $\xi$. In contrast to the FRSB prediction,
   in the spherical case, there exist models $\xi$ such that  $ \# \{ \supp(\gamma_{\star}) \} = k$, for
    $k=1,2$ \cite{AZeng}.
(Here we continue to assume no magnetic field, i.e. no linear field in the Hamiltonian $H_N$.)

Roughly speaking, the connection between $\gamma_{\star}$ and the geometry of near optima is
in the fact that the support of $\gamma_\star$ consists of all the values taken
by the `overlap' $\<\bsigma^1,\bsigma^2\>/N$, when $\bsigma^1 ,\bsigma^2\in\Lev_N(\eta)$
and $\eta$ is close to $\OPT$. 

In order to formalize this statement, for fixed $\eta >0$ and Borel measurable set $A\subset [-1,1],$ set
	\begin{align}
\rP_N(\eta,A)&:=\prob\bigl(\exists \; \bsigma^1 ,\bsigma^2\in \Lev_N(\eta) \mbox{ with } \<\bsigma^1,\bsigma^2\>/N\in A
\bigr).
	\end{align}	
	In other words, $\rP_N(\eta,A)$ is the probability that there exist two spin configurations 
	with energy values above $N \eta$ and whose overlap lies in $A$.	
	\begin{theorem}[\cite{ArnabChen15, AC18Advances}]\label{thm:GammaStar} 
	 Consider $\Sigma_{N} = \{+ 1,-1 \}^{N}$ or $\Sigma_{N} = \S^{N-1}(\sqrt{N})$,
	 and the Hamiltonian $H_N$, with no magnetic field. Assume $\xi$ to be even. Let $u\in[-1,1]$ with $|u|\in\supp(\gamma_{\star})$. 
	 Then, for any $\varepsilon>0,\eta < \OPT $, there exists $K>0$ such that for all $N\geq 1$,
	 
				\begin{align}\label{thm1:eq1}
				\rP_N\bigl(\eta,(u-\varepsilon,u+\varepsilon)\bigr)&\geq 1-Ke^{-\frac{N}{K}}.
				\end{align}
	\end{theorem}

One should expect that for values of $u$ outside the support of $\gamma_{\star}$ the opposite holds. 
Namely, for $\eta$ sufficiently large but smaller than $\OPT$, the probability of finding 
two configurations
 $\bsigma^1,\bsigma^2\in\Lev_N(\eta)$ with overlap $u$ such that $|u|\not\in\supp(\gamma_{\star})$
 should approach zero as $N\to\infty$. This is only known in the case of either 
 positive external field or in the spherical model under the extra assumption that the 
 $1$-RSB occurs at zero temperature\footnote{Namely, $\gamma_{\star}(t)= c_0+c_1\bfone(t\ge t_0)$
 for some constants $c_0,c_1,t_0$.} \cite{AC18Advances}. The last condition is known to hold
 (in the spherical model)
  for instance if the model satisfies 
				$\xi'(1)>\xi''(0)(1+z)$ and $s/\xi'(s)$ is convex on $(0,1)$, where $z$ is the unique solution of
	\begin{align}
	\frac{1}{\xi'(1)}&=\frac{1+z}{z^2}\log (1+z)-\frac{1}{z}.
	\end{align} 

We conclude this section by mentioning that both in the spherical and in the Ising model, it is known 
that there exist exponentially many nearly orthogonal configurations near the ground state energy 
\cite{ArnabChen15, AC18Advances}. 
More precisely, for $\epsilon,\eta,K>0$ and $q\in[0,1],$ denote by 
   $ \overline{\rP}_{N}(\eps,\eta,q,K)$
    the probability that there exists a subset $O_N\subseteq \Sigma_{N}$ such that
    \begin{enumerate}
    	\item[$(i)$] $O_N\subset \Lev_N(\eta)$.
    	\item[$(ii)$] $O_N$ contains at least $Ke^{N/K}$ many elements.
    	\item[$(iii)$] $|\<\bsigma,\bsigma'\>|\leq N\eps$ for all distinct $\bsigma,\bsigma'\in O_N.$
    \end{enumerate}
Then is its known that for any $\eps,\eta>0$, there exists $K>0$ such that for any $N\geq 1,$
\[
		\overline{\rP}_{N}(\varepsilon,\eta,q_0,K)\geq 1-Ke^{-N/K}.
\]
This phenomenon is connected to chaos in disorder and  was explored by many authors. 
The interested reader can consult \cite{chattbook, DEZ, AC18Advances}.

\section{The TAP approach  and free energy}
\label{sec:Free}

Can we approximate the superlevel sets $\Lev_N(\eta)$ by simpler subsets
of $\Sigma_N$? One way to make this question more precise is to consider an explicit family
of subsets of $\Sigma_N$. It turns out that a convenient such family is provided by
thin  `bands' of the form
\[
\B(\bm,\delta):=\big\{
\bsigma\in\Sigma_N: |\<\bsigma-\bm,\bm\>|<N\delta
\big\},
\] 
around points $\bm\in\oSigma_N:=\conv(\Sigma_N)$ with small $\delta>0$. We then can ask 
how can we choose $\bm\in \oSigma_N$ so that the volume of 
$\B(\bm,\delta)\cap\Lev_N(\eta)$ is approximately the same as the volume of $\Lev_N(\eta)$.
Since the volumes of $\B(\bm,\delta)$ and $\Lev_N(\eta)$ are exponentially smaller than 
the volume of $\Sigma_N$, it makes sense to require that the approximation holds on the 
logarithmic scale:
\begin{align}
\log\Vol(\B(\bm,\delta)\cap\Lev_N(\eta))= \log\Vol(\Lev_N(\eta))+o(N)\, .\label{eq:Volume}
\end{align}
Note that, if we could answer this question for all $\eta<\OPT$, we would also be
able to identify near optima of $H_N$.

Note that Eq. \eqref{eq:Volume} requires that $\B(\bm,\delta)$ has
measure not exponentially small under the uniform measure on $\Lev_N(\eta)$.
An essentially equivalent way to formulate the same condition is to consider the
Gibbs measure of Eq.~\eqref{eq:Gibbs} and ask that
\begin{align}
\log G_{N,\beta}(\B(\bm,\delta)) = o(N)\, .
\end{align}
This can be restated in terms of the free energy of Eq.~\eqref{eq:FirstFree}.
The question is then for which points $\bm\in\oSigma_N$, we have
\begin{align}
F_{N,\beta} &	= F_{N,\beta}(\bm,\delta) +o_N(1)\, , \label{eq:Gband}\\
F_{N,\beta}(\bm,\delta)& :=
\frac1N\log\Big\{\int_{\B(\bm,\delta)}e^{\beta H_N(\bsigma)}\de\bsigma\Big\}\, .\nonumber
\end{align}

After giving it some thought, one realizes that the set of such points $\bm$ can be very large, 
and it might be a good idea to look for a subset of `special' such points.  
The basic problem is that if there is a small region of $\Sigma_N$ containing a large portion 
of the Gibbs measure $G_{N,\beta}$, there 
are many possible ways to choose $\bm$ such that the band $\B(\bm,\delta)$
will contain most of this region. As an extreme example, consider
the case in which $G_{N,\beta}$ concentrates at a single configuration $\bsigma_0$.
Then, all vectors $\bm\in\oSigma_N$ such that $\|\bm-\bsigma_0/2\|^2_2 \approx \|\bsigma_0/2\|_2^2$
would be selected.

In light of this example, we could require that the Gibbs measure on the band $\B(\bm,\delta)$ 
does not concentrate in any specific direction within the band. Namely, that for any non-random 
$\bsigma'\in \B(\bm,\delta)$ and small $\eps>0$,  
\begin{equation*}
	G_{N,\beta}\Big(
	\,\big|\<\bsigma, \bsigma'\>-\<\bm,\bm\>\big|<\eps N \,\Big|\,\bsigma\in \B(\bm,\delta)
	\Big) = 1- o_N(1)\, .
\end{equation*}
For reasons that will become clearer below, we will work  with a weaker condition.
We require that given two constants 
 $\eps>0$ (small) and $k$ (large), 
\begin{equation}
	\label{eq:orthpts}
	\frac1N\log G_{N,\beta}^{\otimes k}\Big(\forall i\neq j:\,
	\,\big|\<\bsigma^i,\bsigma^j\>-\<\bm,\bm\>\big|<\eps N\,\Big|\,\bsigma^i\in \B(\bm,\delta)
	\forall i\Big) = o_N(1)\, .
\end{equation}
Defining $\tilde\bsigma:=\bsigma-\bm$, for small $\delta$ and $\eps$ the above means 
that the probability to sample many replicas from $\B(\bm,\delta)$ such that 
$\<\tilde\bsigma^i,\tilde\bsigma^j\>/N =o(1)$ for all $i\neq j$ is not exponentially small.

With this extra constraint, can we identify the heavy bands which satisfy both 
Eq.~\eqref{eq:Gband} and Eq.~\eqref{eq:orthpts}? It turns out that the answer is positive, and that 
they are characterized asymptotically by the simple approximate equality
\begin{equation}
	\label{eq:char}
\frac{\beta }{N}H_N(\bm)+\TAP_\beta(\mu_{\bm}) = F_{N,\beta} +o_N(1),
\end{equation}
where $\TAP_\beta(\mu_{\bm})$ is a deterministic functional of the empirical measure 
$\mu_{\bm}:=N^{-1}\sum_{i} \delta_{m_i}$.

Recall that we started this section by asking whether we can approximate the level sets
$\Lev_N(\eta)$ (or the Gibbs measure $G_{N,\beta}$) by simpler subsets of $\Sigma_N$: thin bands
$\B(\bm,\delta)$.
At first sight, this question seems significantly harder than the one of finding near maxima
of $H_N(\bsigma)$. However, Eq.~\eqref{eq:char} suggests that we can 
reduce the problem of finding the optimal band centers $\bm$ to the one
of maximizing the free energy functional $\beta H_N(\bm)+N\TAP_\beta(\mu_{\bm})$.

As implied by the notation, the functional $\TAP_\beta(\mu_{\bm})$, which was introduced
 in \cite{subag2018landscape,chen2018generalized}, is a generalization of the 
 TAP free energy invented by Thouless, Anderson and Palmer \cite{thouless1977solution}. We 
 will define it and explain how it is computed in the Section \ref{sec:TAPcorrection}.
By maximizing over $\bm$, the characterization \eqref{eq:char} above also leads to a 
generalized TAP representation for the free energy which we will state in  
Section \ref{sec:TAPrepresentation}. Finally, in Section \ref{sec:tree} we will describe 
how one can obtain from the classical pure states decomposition a tree all of whose 
vertices satisfy \eqref{eq:char}, which will motivate some of the algorithms in 
Section \ref{sec:algorithms}.

\subsection{Generalized TAP correction} \label{sec:TAPcorrection}

Define the set of $k$-tuples from the band $\B(\bm,\delta)$:
\[
\B_k(\bm,\delta,\eps):=\Big\{
(\bsigma^1,\ldots,\bsigma^k)\in \B(\bm,\delta)^k: \forall i\neq j,\, 
|\<\bsigma^i,\bsigma^j\>-\<\bm,\bm\>|<N\eps
\Big\}
\]
and the replicated free energy
\begin{equation}
	\label{eq:TAPk}
	\TAP_{N,\beta,k}(\bm,\delta,\eps):=\frac{1}{kN}
	\log\Big\{\int_{\B_k(\bm,\delta,\eps)}e^{\beta\sum_{i\leq k}[H_N(\bsigma^i)-H_N(\bm)] }
	\de\bsigma^1\cdots \de\bsigma^k\Big\}.
\end{equation}
Observe that
\begin{equation}
	\label{eq:TAPineq}
	F_{N,\beta}(\bm,\delta)=\frac{\beta }{N}H_N(\bm)+\TAP_{N,\beta,1}(\bm,\delta,\eps)\geq 
	\frac{\beta }{N}H_N(\bm)+\TAP_{N,\beta,k}(\bm,\delta,\eps),
\end{equation}
where the difference of the two sides of the inequality is exactly the left-hand side of
 \eqref{eq:orthpts} divided by $k$.

Define the limit
\begin{align}
\label{eq:TAP-Def-Lim}
\TAP_\beta(\mu):=\inf_{\delta,\eps,k} \lim_{N\to\infty}\E \, \TAP_{N,\beta,k}(\bm_N,\delta,\eps),
\end{align}
where $\bm_N$ is an arbitrary sequence such that $\mu_{\bm_N}\Rightarrow\mu$. 
The fact that the limit in Eq.~\eqref{eq:TAP-Def-Lim} exists is of course highly 
non-trivial and is established in \cite[Theorem 1]{chen2018generalized} for Ising models and in \cite[Propostion 1]{subag2018landscape} for spherical models.\footnote{To be precise, what is proved is that with the limit in $N$ 
replaced either by $\limsup$ or $\liminf$ the above converges to the same quantity. 
Technically, in the definition above we should use either of them.} Below we 
will give more explicit expressions for $\TAP_\beta(\mu)$ and sketch how 
it is computed. Before that, we discuss one of its key properties:
the fact that it concentrates uniformly in $\bm$.

For a single point $\bm$ and any $\delta$, $\eps$ and $k$, the deviation of
 $\TAP_{N,\beta,k}(\bm,\delta,\eps)$  from its mean is small with
   high probability, as for the usual free energy. More precisely,
   standard concentration-of-measure arguments imply that, 
   for any $\bm\in\oSigma_N$ and any $k,\eps,\delta,t>0$ there exists $c>0$
   such that, for all $N$ large enough
   \begin{align}
   	\P\Big(\big|\,
	\TAP_{N,\beta,k}(\bm,\delta,\eps) - \E \TAP_{N,\beta,k}(\bm,\delta,\eps)
	\,\big|<t
	\Big)>1-e^{-cN}.
	\end{align}
   The maximal deviation over all
    $\bm$, however,  is typically of order $O(1)$.  
    
    However, the concentration of $\TAP_{N,\beta,k}(\bm,\delta,\eps)$ improves as $k$
    gets larger and $\eps,\delta$ get smaller.
Indeed, the strength of the concentration is dictated by the  maximal variance of the process  
$\sum_{i\leq k}[H_N(\bsigma^i)-H_N(\bm)]$ over  $\B_k(\bm,\delta,\eps)$ divided by
 $k^2$, (coming from the division by $k$ before the $\log$) which is of the order of 
 $N(1/k+\delta+\eps)$. As we let $\delta,\eps\to0$ and $k\to\infty$ we have the
   following uniform concentration result. 
\begin{theorem}[Generalized TAP correction {\cite[Theorem 1]{chen2018generalized}}, {\cite[Proposition 1]{subag2018landscape}}] \label{thm:concentration}For any $c,t>0$, if $\delta,\eps>0$ are small enough and $k\geq1$ is large enough, then for large $N$,
	\[
	\P\Big(
	\forall \bm\in \Chull:\,\big|\,
	\TAP_{N,\beta,k}(\bm,\delta,\eps) - \TAP_\beta(\mu_{\bm})
	\,\big|<t
	\Big)>1-e^{-cN}.
	\]
\end{theorem}

Next we discuss the computation of  $\TAP_\beta(\mu)$. Recall the notation 
$\tilde\bsigma=\bsigma-\bm$. If $\bsigma\in \B(\bm,\delta)$, for small $\delta$, 
\begin{equation}\label{eq:1spin}
	H_N(\bsigma)-H_N(\bm) = H_N^{\bm}(\tilde\bsigma) + \nabla H_N(\bm)\cdot \tilde\bsigma +NO(\delta)\, .
\end{equation}
We can think of  $H_N^{\bm}(\tilde\bsigma)$ as the residual of the first order Taylor
expansion of $H_N(\bsigma)$ around $\bm$. 
This turns out to be a Gaussian spin glass Hamiltonian, with covariance
\begin{align}
\E H_N^{\bm}(\tilde\bsigma^1)H_N^{\bm}(\tilde\bsigma^2)=N \tilde\xi_q(\<\tilde\bsigma^1,\tilde\bsigma^2\>/N)\, .
\end{align}
Here $\tilde\xi_q(x):=\xi(x+q)-\xi(q)-\xi'(q)x$. Note that
$\xi_q(0)=\xi_q'(0)=0$ corresponding to the fact
that the Hamiltonian  $H_N^{\bm}(\tilde\bsigma)$ which does not contain 
linear terms in $\tilde\bsigma$ (a `random magnetic fiels' or `$1$-spin interaction' in physics language). 

Denote by $\TAP^0_{N,\beta,k}(m,\delta,\eps)$ the free energy defined as in Eq.~\eqref{eq:TAPk},
 but with $H_N(\bsigma^i)-H_N(\bm)$ replaced by $H_N^{\bm}(\tilde\bsigma^i)$, for each $i\leq k$. 
 Note that if $|\<\tilde\bsigma^i,\tilde\bsigma^j\>|<N\eps$ for any $i<j\leq k$, then
\begin{equation}
	\label{eq:extbd}
	\Big|\frac1k\sum_{i\leq k} \nabla H_N(\bm)\cdot \tilde\bsigma^i\Big| 
	\leq 
	\Big\|\frac1k \sum_{i\leq k} \tilde\bsigma^i \Big\| \cdot\big\|\nabla H_N(\bm)\big\|\leq \sqrt{N(1/k+\eps)}\big\|\nabla H_N(m)\big\|=\sqrt{1/k+\eps}\cdot O(N).
\end{equation}
Hence,  
\begin{equation}
	\label{eq:TAPext}
	\E \,\TAP_{N,\beta,k}(m,\delta,\eps) = \E \,\TAP^0_{N,\beta,k}(\bm,\delta,\eps)
	+o_{\delta,\eps,k}(1)\, ,
\end{equation}
where $o_{\delta,\eps,k}(1)$ denotes an error bounded uniformly in $N$ by a quantity that vanishes if we
let $\delta,\eps\to 0$ and $k\to\infty$. 

At this point, the treatment is slightly different in the spherical and Ising cases.

\emph{In the spherical case $\Sigma_N=\S^{N-1}(\sqrt{N})$,} it can be shown that for small 
$\delta$ and any $\eps$ and $k$,
\begin{equation}\label{eq:KeyTAP_Spherical}
	\E \,\TAP^0_{N,\beta,k}(\bm,\delta,\eps) = \E\, \TAP^0_{N,\beta,1}(\bm,\delta) +o_{\delta}(1),
\end{equation}
where $o_{\delta}(1)$ denotes an error bounded uniformly in $N$ by a quantity 
independent of $k,\eps$, that vanishes if we let $\delta\to 0$.
Note that $\TAP^0_{N,\beta,1}(\bm,\delta):=\TAP^0_{N,\beta,1}(\bm,\delta,\eps)$ does not depend 
on $\eps$, as it is computed with one replica.
That is, once we remove the external field in \eqref{eq:1spin}, computing the free energy 
with many replicas instead of one has no cost asymptotically.

By combining the above one obtains that the TAP correction only depends on
$\mu$ via its second moment 
$q=\int \! x^2\, \mu(\de x)$ and, for an arbitrary sequence $\bm=\bm_N$ such that 
$\<\bm_N,\bm_N\>/N\to q$, is given by
\begin{align}
\TAP_\beta(\mu)=\inf_{\delta}\lim_{N\to\infty}\E \TAP^0_{N,\beta,1}(\bm,\delta)=
\frac12\log(1-q)+F_\beta(q)\, .\label{eq:TAP-Spherical}
\end{align}
The logarithmic term, which accounts for a change of volume, is equal to the 
asymptotic entropy 
\begin{align}
\inf_{\delta}\lim_{N\to\infty}\frac{1}{N}\log \Vol(\B(\bm_N,\delta)) =\frac12\log(1-q)
\end{align}
The term  $F_\beta(q)$ is more interesting and is defined as the free energy of the spherical spin glass model with mixture 
$\xi_q(x):=\tilde\xi_q((1-q)x)$. 
The expression of Eq.~\eqref{eq:TAP-Spherical} simplifies when $q= q_{\sEA}$ 
is the Edwards-Anderson parameter, i.e. the rightmost point in the support of 
the Parisi measure. In this case, the model for $\tilde\bsigma$ is in its
paramagnetic phase, which means that its partition function
is well approximated by its expectation: in particular, this submodel is replica symmetric.
Hence, we get the simple
expression  $F_\beta(q_{\sEA})=\frac12\beta^2\xi_{q_{\sEA}}(1)$,
 see \cite[Corollaries 5,11]{subag2018landscape}, which is known as 
 the Onsager correction term and appears in the original TAP formulation.

\emph{In the Ising case $\Sigma_N=\{+1,-1\}^N$,} computing 
$\E \,\TAP^0_{N,\beta,k}(\bm,\delta,\eps)$ is more challenging,
as Eq.~\eqref{eq:KeyTAP_Spherical} does not hold anymore.
In other words, removing the `magnetic field' term is 
not enough to make the replicas roughly independent.

In this case, for a function $h:[-1,1]\to\R$, denote by $\TAP^h_{N,\beta,k}(\bm,\delta,\eps)$ 
the free energy computed as in Eq.~\eqref{eq:TAPk}, but with $H_N(\bsigma^i)-H_N(\bm)$
 now replaced by $H_N^{\bm}(\tilde\bsigma^i)+\sum_{j\leq N} h(m_j)\tilde\bsigma^i_j$. By a 
 bound similar to \eqref{eq:extbd}, the approximation as in \eqref{eq:TAPext} also holds for 
 such general $h$ instead of $h\equiv 0$. Using Eq.~\eqref{eq:TAPineq},
  we therefore have that for large $k$,
\[
\E \,\TAP_{N,\beta,k}(\bm,\delta,\eps) = \E \,\TAP^h_{N,\beta,k}(m,\delta,\eps)  +
o_{\delta,\eps,k}(1)
\leq \E\, \TAP^h_{N,\beta,1}(\bm,\delta)+o_{\delta,\eps,k}(1)\, .
\]
It is possible to show that there exists a specific choice of an external field $h$ such
 that the last inequality asymptotically holds as an equality, and therefore
\[
\TAP_\beta(\mu)=\inf_{h,\delta} \lim_{N\to\infty}\E\, \TAP^h_{N,\beta,1}(\bm_N,\delta),
\] 
where $\bm_N\in[-1,1]$ is an arbitrary sequence such that $\mu_{\bm_N}\Rightarrow\mu$. 

The last result was proved in \cite{chen2018generalized}, which also prove that the TAP correction
admits the following  Parisi type formula,
\[
\TAP_\beta(\mu)=\inf_{\zeta\in\cuP([0,1])}
\Big(
\int \Lambda_{\zeta}(q,a)\, \mu(\de a)-\frac12\beta^2\int_q^1s\xi''(s)\zeta(s)\, \de s
\Big),
\]
where $q=\int\! x^2\, \mu(\de x)$, $\cuP([0,1])$ is the space of probability measures on $[0,1]$,
$
\Lambda_\zeta(q,a):=\inf_{x\in\R}\big(
\Phi_\zeta(q,x)-ax
\big)
$ 
and $\Phi_\zeta(q,x)$ is the solution of the Parisi PDE
\[
\partial_t\Phi_\zeta=-\frac12\beta^2\xi''(t)\Big(
\partial_{xx}\Phi_\zeta+\zeta([0,t])(\partial_x\Phi_\zeta)^2
\Big)
\]
on $[0,1]\times\R$ with boundary condition $\Phi_\zeta(1,x)=\log2\cosh x$.
This is the non-zero temperature version of the Parisi PDE of Eq.~\eqref{eq:PDEFirst}.

Again, the expression for $\TAP_\beta(\mu)$ simplifies when $q=q_{\sEA}$, in certain cases.
Namely, under a certain optimality condition which implies Plefka's condition 
\cite{plefka1982convergence},  which holds at least for some points $\bm$ with $\|\bm\|^2/N=q=q_{\sEA}$,
it was proved in \cite{chen2018generalized} that the correction coincides with 
the classical Onsager correction 
\[
\TAP_\beta(\mu)=
-\int_{\reals}\! I(x)\,\mu(\de x)+\frac12\beta^2\xi_q(1)\, .
\]
Here $I(x)=\frac{1+x}{2}\log\frac{1+x}{2}+\frac{1-x}{2}\log\frac{1-x}{2}$ and the
 first term above corresponds to the entropy of the band as before.

\subsection{Generalized TAP representation}
\label{sec:TAPrepresentation}

For $q$ in the support of the Parisi measure, there exists a (random) sequence
 $\bm_\star=\bm_{\star,N}$ with $\<\bm_\star,\bm_\star\>/N=q$ such that for any 
 $\eps$, $\delta$ and $k$, with high probability  both \eqref{eq:Gband} and \eqref{eq:orthpts}
  hold, or equivalently,
\begin{equation}
	\label{eq:TAPapx}
	F_{N,\beta} = F_{N,\beta}(\bm_\star,\delta) +o_N(1) =
	 \frac{\beta }{N}H_N(\bm_\star)+\TAP_{N,\beta,k}(\bm_\star,\delta,\eps) +o_N(1)
	=\frac{\beta }{N}H_N(\bm_\star)+\TAP_{\beta}(\mu_{\bm_\star}) +o_N(1),
\end{equation}
where  the last approximation follows from the concentration of Theorem \ref{thm:concentration}. 
By \eqref{eq:TAPineq}, for any $m$,
\begin{equation*}
	F_{N,\beta} \geq \frac{\beta }{N}H_N(\bm)+\TAP_{\beta}(\mu_{\bm})+o_N(1).
\end{equation*}
The following theorem follows.

\begin{theorem}[Generalized TAP representation {\cite[Theorem 2]{chen2018generalized}}, {\cite[Theorem 4]{subag2018landscape}}] \label{thm:TAPrepresentation}
	For any $q$ in the support of the Parisi measure, in probability,
	\begin{equation}
	\label{eq:TAPrep}
	\lim_{N\to\infty}\Big|
	F_{N,\beta} - \max_{\|\bm\|_2^2=Nq}\Big(\frac{\beta}{N}H_N(\bm)+\TAP_\beta(\mu_{\bm})\Big)
	\Big|=0.
	\end{equation}
\end{theorem}
An explicit expression using the Parisi measure for the value of $\frac1NH_N(\bm)$
 at any approximate maximizer of the above can  also be derived, see 
  \cite[Proposition 10]{subag2018landscape} and \cite[Theorem 5]{chen2020generalized2}.

In the spherical case, $\TAP_\beta(\mu_{\bm})$ is constant on the sphere and only the Hamiltonian
 is maximized in \eqref{eq:TAPrep}. For the pure spherical models, this can be used to compute 
 the free energy at any temperature \cite{subag2021pureTAP}, also for multi-species models 
 \cite{subag2021multi1,subag2021multi2}. For earlier results on the TAP representation 
 with the classical Onsager correction see 
 \cite{belius2019TAPplefka,benarous2021shattering,benarous2020geometry,bolthausen2019morita,
 chen2018TAP,crisanti1995TAP,kurchan1993barriers,subag2017geometry}.

\begin{remark}
It is useful to comment on the implication of Theorem \ref{thm:TAPrepresentation}
for optimization. In this case we are interested in maximizing $H_N(\bm)$ subject to
$\bm\in \oSigma_N$ and $\<\bm,\bm\>/N =1$ (the last two constraints are equivalent to 
$\bm\in\Sigma_N$).

Since we are only interested in achieving any value $\eta<\OPT$, it is reasonable to
perturb this problem in two ways. First, we will accept $\<\bm,\bm\>/N \ge q_\star$
for some $q_\star$ close to $1$. Second, we perturb the cost function $H_N(\bm)$
by considering 
$H_N(\bm)+(N/\beta)\TAP_\beta(\mu_{\bm})$.

Theorem \ref{thm:TAPrepresentation} states a remarkable property of this modified objective
function. Namely, for any $q$ in the support of Parisi's measure (and therefore, possibly, for
any $q\in[0,q_{\star}]$) this objective function has an optimizer on the sphere 
of radius $\sqrt{Nq}$, with value close to the optimum on the sphere $\sqrt{Nq_\star}$.
This naturally suggests to find the latter by gradually increasing this radius.
This  is indeed what we will do in Section \ref{sec:algorithms}, although not explicitly using 
the free energy $H_N(\bm)+(N/\beta)\TAP_\beta(\mu_{\bm})$.
\end{remark}

\subsection{The tree of pure states}
\label{sec:tree}

Recall the pure state decomposition of the Gibbs measure
 \cite{mezard1987spin,talagrand2010states,jagannath2017apxult}.
 Namely, at low temperatures there exists (as proven in \cite{talagrand2010states,jagannath2017apxult} for
 `generic' mixtures $\xi$) 
 a partition $\Sigma_N=\cup_{\alpha\in \cS} B_{\alpha}\cup D$,
 so that $G_{N,\beta}(D) = o_N(1)$  ($D$ is negligible) and therefore
\begin{equation}\label{eq:PureStatesDecomposition}
	G_{N,\beta}(\,\cdot\,) \approx \sum_{\alpha\in \cS} w_\alpha\, G^{(\alpha)}(\,\cdot\,)\, ,
	\;\;\;\;\; w_\alpha:= G_{N,\beta}(B_\alpha)\,  ,\;\;
	G_{N,\beta}^{(\alpha)}(\,\cdot\,) := G_{N,\beta}(\,\cdot\, |\,B_\alpha)\, .
\end{equation}
Further, each `pure state' $G_{N,\beta}^{(\alpha)}$ 
can be characterized by its barycenter (magnetization)
\[
\bm^{(\alpha)}:=\int\!\bx \,G_{N,\beta}^{(\alpha)}(\de\bx)\, .
\]
Namely,  uniformly in $\alpha,\alpha'\in\cS$, for any $\eps>0$,
\begin{equation}
	\label{eq:pure}
	G^{(\alpha)}_{N,\beta}\otimes G^{(\alpha')}_{N,\beta}
	\big(|\<\bsigma^1,\bsigma^2\>-\<\bm^{(\alpha)},\bm^{(\alpha')}\>|\ge N\eps\big)
	= o_N(1)\, .
\end{equation}
If the Parisi distribution has an atom at $q_{\sEA}$ then the weights $w_\alpha$ are of order 
$O(1)$, otherwise they decay to zero sub-exponentially as $N\to\infty$. 
In both cases, their asymptotic distribution can be described by  Ruellle probability 
cascades or a limit of such.

The pure states are arranged ultrametrically 
\cite{mezard1984nature1,mezard1984nature2,panchenko2013parisi}:
namely for any three pure states $\alpha_1,\alpha_2,\alpha_3\in\cS$,
letting $R_{ij}=N^{-1}|\<\bm^{(\alpha_i)},\bm^{(\alpha_j)}\>|$,  
we have $R_{12} = \max(R_{13},R_{23})+o_N(1)$, where the term $o_N(1)$ is uniform 
in $\alpha_{1},\alpha_2,\alpha_3$. 
As a consequence, we can cluster hierarchically the pure states, and associate
 to each cluster
  $C$ the average $\bm^{(C)}:=\frac{1}{|C|}\sum_{\alpha\in C}\bm^{(\alpha)}$.
  
This hierarchical structure can be represented
 by a rooted tree $\cT_{N,\beta}$.
 Leaves of the tree are indexed by $\alpha \in \cS$ and associated to pure states 
 centers $\bm^{(\alpha)}$. 
 Each other vertex of the tree corresponds to a cluster. Clusters will be indexed by 
 $\alpha\in \cV\setminus \cS$,  and are associated to the cluster center $\bm^{(\alpha)}\in\Chull$.
 The theorem below summarizes a few properties of this structure.
 We denote by $\alpha_1 \wedge \alpha_2$ the least common ancestor of $\alpha_1$ and $\alpha_2$, by 
 $|\alpha|$ 
the depth of $|\alpha|$ (the depth of the root being $0$, and write $\alpha_1\preceq \alpha_2$ if 
$\alpha_2$ is a descendant of $\alpha_1$ (including $\alpha_1=\alpha_2$). 
\begin{theorem}[Ultrametric tree {\cite[Corollary 13]{subag2018landscape}}] \label{thm:Tree}
If $\xi(t)=\sum_{k\ge 2}c_k^2t^k$ is generic\footnote{Namely, there are infinitely many values of 
$k$ even such that $c_k\neq 0$ and infinitely many values of
$k$ odd such that $c_k\neq 0$.}
 and  $\beta$ is in the RSB phase\footnote{Namely, Parisi's measure is not 
 a single point mass.}, then there exist a decomposition of the Gibbs
 measure of the form \eqref{eq:PureStatesDecomposition} and an associated rooted
 tree $\cT_{N,\beta}$ with vertex set $\cV=\cV_{N,\beta}\in \Chull$, 
 and associated vectors $\{\bm^{(\alpha)}\}_{\alpha\in\cV}$ with the following properties:
  \begin{enumerate}
 \item  The leaves of $\cT_{N,\beta}$ are the barycenters $\{\bm^{(\alpha)}:\,\alpha\in \cS\}$ 
 of the pure states $G^{(\alpha)}_{N,\beta}$.
\item   The tree is regular and complete (i.e. all leaves have the same depth)
and its degree diverges.
\item  For any $\alpha\in \cV$, the value of $\|\bm^{(\alpha)}\|_2^2/N$ only depends on the depth 
$|\alpha|$.
\item  In the $k$-RSB case (that is, if Parisi's measure consists of $k+1$
point masses) the depth of $T$ is equal to $k$.  Further $Q:=\{\|\bm^{(\alpha)}\|_2^2/N:\,\alpha\in \cV\}$ is equal 
    to the support of the Parisi measure.
\item In the full-RSB case (that is, if the support Parisi's measure has infinite cardinality) 
the depth of the tree diverges and $Q$ converges to the support of Parisi's measure 
as $N\to\infty$.
\item 
Finally, the following hold
with probability going to $1$, for all $\alpha_1,\alpha_2\in \cV$ and 
$\alpha\in \cS$ (for a deterministic $\delta_N=o_N(1)$):
\begin{enumerate}
	\item $\displaystyle  \bm^{(\alpha_1)}\preceq \bm^{(\alpha)}\implies B_\alpha\subset  \B(\bm^{(\alpha_1)},\delta_N).$
	\item $\displaystyle \<\bm^{(\alpha_1)},\bm^{(\alpha_2)}\>/N =\|\bm^{(\alpha_1\wedge\alpha_2)}\|^2_2/N +o_N(1)).$
	
	\item \label{item:3}$\displaystyle 
	\beta N^{-1}H_N(\bm^{(\alpha_1)})+\TAP_\beta(\mu_{\bm^{(\alpha_1)}}) = F_{N,\beta}+o_N(1)$.
\end{enumerate}  
\end{enumerate}
\end{theorem}

 \begin{remark}
Except for the last part involving the TAP correction, this theorem follows from the decomposition
 of \cite{jagannath2017apxult}.
  For inner vertices, the last point follows 
 since the corresponding band contains many orthogonal pure states.
 \end{remark}

\section{Algorithms}
\label{sec:algorithms}

\subsection{Typical case approximation}

Is there a polynomial-time algorithm that finds near optima of the 
spin glass Hamiltonian $H_N(\bsigma)$?

From a worst case perspective, the answer is likely to be negative. 
Consider, to be definite, the pure degree-$k$
Hamiltonian  $H^{(k)}_N(\bsigma) = N^{-(k-1)/2}\<\bG^{(k)},\bsigma^{\otimes k}\>$.
In theoretical computer science, the problem of maximizing this function over 
$\Sigma_N$ is formulated in terms of existence of an `approximation algorithm.'
An approximation algorithm  accepts as input $\bG^{(k)}$ and outputs $\bsigma^{\salg}\in
\Sigma_N$ such that, for any input $\bG^{(k)}$, we are guaranteed  to have
\begin{align}
 \max_{\bsigma\in\Sigma_N} H^{(k)}_N(\bsigma) \ge 
 H^{(k)}_N(\bsigma^{\salg})\ge \rho \cdot \max_{\bsigma\in\Sigma_N} H^{(k)}_N(\bsigma)\, ,
\end{align}
with $\rho$ independent of $H_N$ (but possibly dependent on $N$). 
It is immediate to show that this can be achieved for the quadratic Hamiltonian
$H^{(k=2)}_N(\bsigma)$ under the spherical constraint $\Sigma_N = \S^{N-1}(\sqrt{N})$.
Indeed in this case the problem can be reduced to an eigenvalue problem that is solved efficiently
for any fixed approximation ratio $\rho\in(0,1)$. 

However, for all the other cases this objective is all but hopeless for any constant $\rho\in (0,1)$.
Indeed, achieving $\rho>1/(\log N)^c$ (for $c$ a small constant) 
is NP-hard for $k=2$ and $\Sigma_{N}=\{+1,-1\}^N$ \cite{arora2005non} 
(the SK model). 
For higher-order polynomials, the task is even more difficult. For instance,  \cite{barak2012hypercontractivity} 
proves that obtaining $\rho> \exp(-(\log N)^c\}$ is NP-hard already for the spherical case
$\Sigma_N = \S^{N-1}(\sqrt{N})$.

These hardness results motivate the search for algorithm that achieve approximation factor
$\rho$ \emph{with high probability} with respect to the realization of the Hamiltonian.
Insights from spin glass theory can be brought to bear on the optimization
question when considering random realizations of the objective.

It is convenient to come back the general case of
a mixed spin glass $H_N(\bsigma)$, as defined in Eq.~\eqref{eq:FirstHamiltonian}.
Formally, we are interested in achieving the following guarantee, in polynomial time:
\begin{align}
  \lim_{N\to\infty}\prob\Big(H_N(\bsigma^{\salg})\ge \rho\cdot
  \max_{\bsigma\in\Sigma_N}H_N(\bsigma)\Big) =1\, .
  \label{eq:Prob}
  \end{align}
The supremum value of $\rho$ such that this is possible will be referred
to as the `typical case' approximation ratio $\rho_*$.
Since  $\max_{\bsigma\in\Sigma_N}H_N(\bsigma)/N$ concentrates 
around its expectation, which has a limit as $N\to\infty$
given by Parisi formula \eqref{eq:ParisiF}, determining such an approximation
ratio (for a specific algorithm), amounts to determine the asymptotics of 
$H_N(\bsigma^{\salg})/N$.

The structure of the space of near optima described in the previous section 
suggests a possible approach to construct approximate ground states.
Summarizing,  near optima  are organized according to a
tree with vertices associated to `magnetization' vectors, i.e. 
vectors $\bm^{(\alpha)}\in \oSigma_N := \conv(\Sigma_N)$.
At zero temperature, leaves of the tree are configurations $\bsigma^{(\alpha)}\in\Sigma_N$,
and near optima are concentrated in the vicinity of such configurations.
As in the rest of the paper, we will assume for simplicity $\xi'(0) =0$
(no magnetic field). In this case the root of the tree can be taken to be the
zero vector $\bzero$.

If the Hamiltonian $H_N(\bsigma)$ presents replica symmetry breaking at zero
temperature with un-normalized Parisi measure $\gamma_\star$, 
then there exist internal nodes of this tree at nearly every radius in the support
of $\gamma_\star$.  Namely, for any $q\in \supp(\gamma_\star)$, there exists
a node $\bm^{(\alpha)}$ of the tree with $\|\bm^{(\alpha)}\|_2^2 =N q+o(n)$.
(See Theorem \ref{thm:Tree}.)

This picture suggests a possibility for constructing an algorithm 
that finds near ground states: if $\supp(\gamma_\star) = [0,1]$, we can hope to
initialize the algorithm with $\hbm^{0} = \bzero$, and then follow a branch of the tree 
by constructing a trajectory $t\mapsto \hbm^t\in \oSigma_N$, indexed by 
the `time' variable $t\in [0,1]$. Of course time will have to be discretized in
the algorithm, and hence we will compute recursively $\hbm^{t+\delta}$
as a function of $\hbm^{\le t}:= (\hbm^s)_{s\le t}$.

One observation turns out to be crucial in the construction of such an algorithm.
If $\alpha'$ is a descendant of $\alpha$ in the tree (See Theorem \ref{thm:Tree}),
then $\bm^{(\alpha')}-\bm^{(\alpha)}$ is roughly orthogonal to $\bm^{(\alpha)}$
(namely, $\<\bm^{(\alpha')}-\bm^{(\alpha)},\bm^{(\alpha)}\>=o(N)$).
It is natural to enforce the same property in the algorithm
updates by requiring that the update $\hbm^{t+\delta}-\hbm^t$ is 
roughly orthogonal $\hbm^t$.

In Sections \ref{sec:HessianDescent} and \ref{sec:IAMP} we will discuss two 
implementations of this idea that yield two algorithms, respectively 
for the spherical case $\Sigma_N=\S^{N-1}(\sqrt{N})$ or the general case
$\Sigma_N\in\{\S^{N-1}(\sqrt{N}),\{+1,-1\}^N\}$.
Before describing these algorithms, we state a result characterizing the energy 
$H_N(\bsigma^{\salg})$ that they achieve (and hence the corresponding approximation ration $\rho$). 
As we will see in Section \ref{sec:Hardness}
this is the optimal energy achieved by any algorithm in a broader class, which includes,
among others,  simulated-annealing type algorithms (when run for physical time independent 
of $N$).

In order to state the characterization of $H_N(\bsigma^{\salg})$,
we introduce the following space of order parameters:
\begin{align}
\cuL := \Big\{\gamma:[0,1)\to \reals_{\ge 0}: \;\; \|\xi''\gamma\|_{\sTV[0,t]}<\infty~ \forall t\in [0,1), \int_0^1\!\xi''\gamma(t)\,\de t < \infty\Big\} \, .\label{eq:LDef}
\end{align}
We endow this space with the weighted $L^1$ metric
\begin{equation}\label{eq:modified_L1}
\|\gamma_1-\gamma_2\|_{1,\xi''} := \|\xi''(\gamma_1-\gamma_2)\|_1 =\int_0^1\xi''(t)|\gamma_1(t)-\gamma_2(t)|\de t,
\end{equation}
hence implicitly identifying $\gamma_1$ and $\gamma_2$ if they coincide for
almost every $t\in [0,1)$.

The usual space of order parameters $\cuU$ appearing in Parisi formula
is given by non-decreasing right-continuous functions $\gamma:[0,1)\to \reals_{\ge 0}$.
Hence $\cuU = \cuL\cap \{\gamma\;\;$ non-decreasing$\;\}$. The asymptotics of the energy density
$H_N(\bsigma^{\salg})/N$ is given by a modified variational principle in which the 
Parisi functional is minimized over the larger space $\cuU$. 
\begin{theorem}[{\cite[Theorem 3]{el2021optimization}}]\label{thm:AlgoVarPrinciple}
Assume that the infimum $\inf_{\gamma\in \cuL} \Par(\gamma)$ is achieved at a function 
$\gamma^\star_{\cuL}\in\cuL$. 
Further denote by $\chi$ the computational complexity of evaluating $\nabla H_N(\bm)$  at a point 
$\bm\in \oSigma_N$,
and by $\chi_1$  the complexity of evaluating one coordinate of $\nabla H_N(\bm)$ at a point 
$\bm\in\oSigma_N$.

Then for every $\eps>0$ there exists an algorithm with complexity at most 
$C(\eps)\cdot (\chi+N)+N\chi_1$ which outputs
$\bsigma^{\salg}\in\Sigma_N$ such that
\begin{align}
\frac{1}{N}H_N(\bsigma^{\salg}) \ge \ALG-\eps\, ,\;\;\;\;\;
\ALG:=\inf_{\gamma\in\cuL}\Par(\gamma)\, ,\label{eq:E-ALG}
\end{align}
with probability converging to one as $N\to\infty$.
\end{theorem}

\begin{remark}\label{rmk:Spherical}
For the spherical model $\Sigma_N=\S^{N-1}(\sqrt{N})$, the Parisi
PDE can be solved explicitly, yielding the following explicit form 
of Parisi's functional  
\begin{align}
\Par(\gamma) & = \inf_{L\ge \int_0^1 \gamma(s)\de s}\frac{1}{2}\int_0^1\left(\xi''(t)\Gamma(t) +\frac{1}{\Gamma(t)}\right)\de t\,,\label{eq:ParisiSpherical}\\
\Gamma(t) & := L-\int_0^t \gamma(s)\de s\,  .
\end{align}
This is minimized at $\Gamma(t) = 1/\sqrt{\xi''(t)}$, yielding the following simple formula 
for the algorithmic threshold:
\begin{align}
\ALG&= \inf \big\{\Par(\gamma) :\; \gamma\in\cuL\,\big\} = 
\int_{0}^1 \sqrt{\xi''(t)}\, \de t\, .\label{SphericalEalg}
\end{align}
\end{remark}

Still focusing on the spherical case, the order parameter $\gamma^{\star}_{\cuL}$ minimizing
$\Par(\gamma)$ over $\cuL$ is given by 
\begin{align}
\gamma^{\star}_{\cuL}(t) = -\frac{\de\phantom{t}}{\de t} \xi''(t)^{-1/2}=
\frac{\xi'''(t)}{2\xi''(t)^{3/2}}\, .
\end{align}
If $t\mapsto \xi'''(t)/\xi''(t)^{3/2}$ is non-decreasing, this obviously minimizes 
$\Par(\gamma)$ over $\cuU$ as well. Therefore, in this case,
there exists a polynomial-time algorithm that achieves a $(1-\eps)$-approximation
of the ground start for every $\eps>0$. 

For the Ising case $\Sigma_N=\{+1,-1\}^N$, the variational principle does not
admit an explicit solution, and we therefore introduce the following condition.
\begin{definition}[No overlap gap at zero temperature]
A mixed $p$-spin model with mixture $\xi$ is said to satisfy the no-overlap gap 
property at zero-temperature
if there exists $\gamma_{\star}\in\cuU$ strictly increasing in $[0,1)$ such that 
 $\Par(\gamma_\star) = \inf_{\gamma\in \cuU} \Par(\gamma)$.

(Equivalently, this condition holds if the unique solution $\gamma_{\star}$ of Parisi's
formula is strictly increasing.)
\end{definition}
\begin{corollary}\label{coro:NoGap}
Assume the no-overlap gap at zero temperature to hold for the mixture $\xi$. 
Then for every $\eps>0$ there exists an algorithm with the same complexity as in 
Theorem \ref{thm:AlgoVarPrinciple} which outputs
$\bsigma^{\salg}\in \Sigma_N$ such that
\begin{align}
\frac{1}{N}H_N(\bsigma^{\salg}) \ge\frac{1}{N}\max_{\bsigma\in\Sigma_N}H_N(\bsigma) - \eps\, ,
\end{align}
with probability converging to one as $N\to\infty$.
\end{corollary}

In particular, it is currently believed that the classical SK model (corresponding to 
$\xi(t)= ct^2$) enjoys the no overlap gap property, and hence a $(1-\eps)$ approximation 
algorithm exist for any $\eps>0$. 

\subsection{Description of the optimization algorithms}

We will next succinctly review the algorithms that achieve the
energy value of  Theorem \ref{thm:AlgoVarPrinciple}, referring to the original papers
 for further detail. 

\subsubsection{The Hessian ascent algorithm}
\label{sec:HessianDescent}

This algorithm was designed \cite{subag2018following} for the spherical case $\Sigma_N = \S^{N-1}(\sqrt{N})$,
and hence operates on a vector $\bm^t\in \oSigma_N=\Ball^N(\sqrt{N})$
(the ball of radius $\sqrt{N}$ in $N$ dimensions). This vector is updated 
at discrete times $t\in\{\delta,2\delta,\dots,1-\delta,1\}$ starting from
the initialization $\bm^\delta=\Unif(\S^{N-1}(\sqrt{N\delta}))$. 
At each $t\in \delta\naturals$, this vector is updated according to:
\begin{align}
\Hess^{\perp}_t & = \bP_{\bm^t}^{\perp} \nabla^2 H_N(\bm^t)\bP_{\bm^t}^{\perp} \, ,\\
\bm^{t+\delta} &= \bm^t+\sqrt{N\delta}\, \bv_1(\Hess^{\perp}_t)\, . \label{eq:HessianIteration}
\end{align}
Here $\nabla^2 H_N$ is the (ordinary,  Euclidean) Hessian of $H_N:\reals^N\to\reals$,
$\bP_{\bm^t}^{\perp} = \id_N-\bm^t(\bm^t)^{\sT}/\|\bm^t\|_2^2$
is the projector orthogonal to the current state $\bm^t$, and $\bv_1(\bM)$
denotes the (unit norm) eigenvector of a matrix $\bM$ corresponding to the largest eigenvalue.
The sign ambiguity of $\bv_1$ can be removed so as to ensure
$\<\bm^{t+\delta}-\bm^t,\nabla H_N(\bm^t)\> \ge 0$ (which is convenient for the proof)
or also at random. 
Any  other degeneracy is removed at random.

By the Pythagorean theorem, $\|\bm^{t+\delta}\|_2^2= \|\bm^{t}\|_2^2+N\delta$,
whence $\|\bm^t\|_2=\sqrt{Nt}$. The algorithm outputs 
$\bsigma^{\salg}:=\bm^1\in \S^{N-1}(\sqrt{N})$.

The rationale for the iteration \eqref{eq:HessianIteration} is relatively easy to understand. 
At each iteration, the algorithm tries to improve the current value by exploiting the curvature of 
the Hessian. However, instead of maximizing the quadratic approximation of the cost at $\bm^t$
along all directions, we only consider directions that are orthogonal to the 
current state $\bm^t$. 

The use of orthogonal updates is inspired by the structure of the space of near optima described 
in the previous sections. 
Technically, it allows to remove dependencies between gains made at different time steps.
Notice indeed that, by rotational invariance, for any $\bm\in\reals^N$, and for any $k$, 
we have (recalling that $H^{(k)}_N(\bsigma) = N^{-(k-1)/2}\<\bG^{(k)},\bsigma^{\otimes k}\>$
is the degree-$k$ component of the Hamiltonian)
\begin{align}
\bP_{\bm}^{\perp}\nabla^2 H^{(k)}_N(\bm)\bP_{\bm}^{\perp}\ed 
\left(\frac{\|\bm\|^2_2}{N}\right)^{(k-2)/2} \sqrt{\frac{k(k-1)}{2N}}(G^{(k)})_{i,j,1,\dots,1}+G^{(k)})_{j,i,1,\dots,1})_{2\le i,j\le N}, .
\end{align}
whence\footnote{We denote by $\bA\sim\GOE(m)$ a matrix from the Gaussian Orthogonal Ensemble,
i.e. a symmetric matrix with $(A_{ij})_{i\le j\le m}$ independent such that
$A_{ii}\sim\normal(0,2/m)$ and $A_{ij}\sim\normal(0,1/m)$ for $i<j$.}, for $\|\bm\|^2_2 = Nt$,
 $\bP_{\bm}^{\perp}\nabla^2 H_N(\bm)\bP_{\bm}^{\perp}\ed \sqrt{\xi''(t)(1-N^{-1})}\, {\sf GOE}(N-1)$.
If $\bm^t$ was independent of $H_N$, the energy increment at each step of the algorithm would
 therefore be
\begin{align}
H_N(\bm^{t+\delta})-H_N(\bm^{t}) & \ge 
\frac{1}{2} N\delta\lambda_1(\bP_{\bm^t}^{\perp}\nabla^2 H_N(\bm^t)\bP_{\bm^t}^{\perp}) +No(\delta)\, ,\\
 \lambda_1(\bP_{\bm^t}^{\perp}\nabla^2 H_N(\bm^t)\bP_{\bm^t}^{\perp}) &\ge 2 \sqrt{\xi''(t)} \delta-o_N(1)
 \, .
 \end{align}
 (The first inequality holds because the contribution of the first order term 
 in the Taylor expansion is non-negative by our choice of the sign of $\bv_1$.
 The second because $\lambda_1({\sf GOE}(N)) = 2-o_N(1)$.)
 A concentration argument justifies treating $\bm^t$ was independent of $H_N$.
 
 By taking $N\to\infty$ first and $\delta\to 0$ afterwards, we obtain 
 Eq.~\eqref{SphericalEalg}. 
 
\subsubsection{Incremental Approximate Message Passing}
\label{sec:IAMP}

A different approach \cite{montanari2019optimization,el2021optimization}
leverages the general analysis of Approximate Message Passing
(AMP) algorithms developed in 
\cite{bolthausen2014iterative,bayati2011dynamics,javanmard2013state,berthier2020state}.
 We begin by describing the general framework
 of \cite{montanari2019optimization} and then describe how it specializes to 
optimization in mean field spin glasses.

 For $t \in \delta\naturals$, let $f_{t} : \reals^{(t/\delta)+1} \to \reals$ be
a real-valued Lipschitz  function, and let $f_{-\delta}  \equiv 0$. For a sequence of vectors
$\z^{\le t}:=(\z^0,\z^{\delta},\cdots,\z^t)$, $\z^s\in \reals^N$ we use the notation
$f_{t}(\z^{\le t})$ for the vector
$(f_{t}(z_i^0,\cdots,z_i^{t}))_{1\le i \le N}$. For a symmetric tensor 
$\bW \in (\reals^N)^{\otimes p}$ and a vector $\u
\in \reals^N$ we denote by $ \bW\{\u\}$ the vector with $i$-th coordinate
$\sum_{1\le i_1,\cdots,i_{p-1}\le N}W_{i,i_1,\cdots,i_{p-1}}u_{i_1}\cdots u_{i_{p-1}}/(p-1)!$

 Let $Z_0 \sim \normal(0,Q_{0,0})$. For each $t\in\delta\naturals$, let
$(Z_{\delta},Z_{2\delta},\cdots,Z_t)$ be a centered Gaussian vector independent of $Z_0$ with covariance
 $Q_{s,t} = \E[Z_sZ_t]$ defined recursively by
\begin{align}\label{eq:cov_amp}
\begin{split}
Q_{s+\delta,t+\delta} &= \xi'\big(\E(M_sM_t)\big) \;\;\;\;\;
M_s:=f_{s}(Z_{\le s}),~~~ s\ge0. 
\end{split}
\end{align}
The message passing algorithm starts with $\z^0$ with coordinates drawn  i.i.d.\ with distribution 
$\normal(0,Q_{0,0})$ independently of everything else.       
The general message passing iteration takes the form
\begin{align}\label{eq:general_amp}
\begin{split}
\z^{t+\delta} &= \sum_{p=2}^{\infty} c_p \bG^{(p)}\{\bm^t\}- \sum_{j=0}^\ell d_{t, s}\bm^{s-\delta}\, ,
\;\;\;\;\; \bm^s =f_s(\z^{\le s})\, ,\\  
d_{t,s} &= \xi''\big( \E\big[M_tM_{s-\delta}\big]\big) \cdot \E\Big[\frac{\partial f_{t}}{\partial z^s}(Z^{\le t})\Big]\, .
\end{split}
\end{align}  
It can be proven that that, for any fixed $t$, the empirical distribution 
$N^{-1}\sum_{i=1}^N\delta_{(z_i^0,\cdots,z_i^t)}$
is asymptotically Gaussian, with the same covariance as $(Z_0,\dots Z_t)$.
This asymptotic normality result, along with the rule \eqref{eq:cov_amp} for evolution of
covariances is known as `state evolution' and was proven (to various degrees of
generality) in
\cite{bolthausen2014iterative,bayati2011dynamics,javanmard2013state,berthier2020state}
for matrices (the case $p=2$), and generalized in \cite{montanari2019optimization}
to tensors.

Hiding some technical conditions, state evolution implies the following.
Given a discrete time stochastic process $(M_t)_{t\in \delta\naturals}$,
that depends causally on the Gaussian process 
$(Z_t)_{t\in \delta\naturals}$ and satisfies Eq.~\eqref{eq:cov_amp}, the AMP construction
provides (efficiently computable) sequences $(\z^t)_{t\in \delta\naturals}$, $(\bm^t)_{t\in \delta\naturals}$, that
are functions of\footnote{In mathematical language,
$(M_t)$ is progressively measurable with respect to $(Z_t)$,
and $(\z^t)$, $(\bm^t)$ are measurable with respect to $H_N$.} 
the disorder $H_N$ and are approximately  isometric to the original processes. Namely
\begin{align}
\frac{1}{N}\<\bm^t,\bm^s\> = \E[M_tM_s] +o_N(1)\, ,\;\;\;\;\; 
\frac{1}{N}\<\bz^t,\bz^s\> = \E[Z_tZ_s] +o_N(1)\, .
\end{align}

We are planning to output $\bsigma^{\salg} = \bm^{t_*}$ for a certain (fixed)
time $t_*$, or $\bsigma^{\salg} = \Round(\bm^{t_*})$  for a suitable rounding procedure
$\Round: \reals^N\to\Sigma_N$. In practice, $\bm^{t_*}$ will lie close to $\Sigma_N$, 
so the rounding operation does not play an important role in the algorithm.

Motivated by the structure of the space of near optima, and the orthogonality 
property of the associated tree, we seek an algorithm with roughy orthogonal increments:
$\<\bm^t-\bm^s,\bm^s\>/N = o(1)$ for all $t>s$.
 This translates into $\E[(M_t-M_s)M_s] =0$. 
 Further, without loss of generality, we can always
assume the normalization $\|\bm^t\|_2^2=N(t+\delta)$, whence $\E[M_t^2]=t+\delta$. 
Using Eq.~\eqref{eq:cov_amp}, this implies that $(Z_t)$ is a Gaussian process 
with covariance $\E[Z_tZ_s] =\xi'(t\wedge s)$. In other words,
$Z_t$ has independent increments, and therefore must be a time change of the Brownian motion.
Indeed we can write it as
 \begin{align}
Z_t = \int_{0}^t\sqrt{\xi''(s)}\, \de B_s\, ,
\end{align}
with $(B_t)_{t\ge 0}$ a standard Brownian motion. We are then left with the task 
of constructing a process $(M_t)$ that causally depend on $(Z_t)$ (equivalently, on $B_t$) 
and satisfies $\E[M_tM_s] = (t\wedge s)_+\delta$ for $t,s\in \delta\naturals$. It turns 
that these constraints nearly\footnote{Indeed, by the martingale representation theorem, the choice is unique 
when we formally set $\delta\to 0$.} fix $M_t$ to be of the form
\begin{align}
 M_t = M_0+\sum_{s\in \{0,\delta,\dots,t-\delta\} }U_s(Z_{s+\delta}-Z_s)\, ,
 \end{align}
 where $M_0\sim \normal(0,\delta)$ and the process $(U_t)$ depends
causally on $(Z_t)$.
 
 We are therefore left with the task of constructing $U_t = F_u(Z_0,\dots,Z_{t})$.
 The asymptotic value achieved by the algorithm can be computed in terms of this process 
 yielding\footnote{The left-hand side converges in probability to a non-random quantity given 
 by the right-hand side.}
\begin{align}
\lim_{N\to\infty}\frac{1}{N}H_N(\bm^t) = \E\Big[\int_0^t\xi''(s) U_s\,\de s \Big] +o_1(\delta)\, .
\end{align} 

The algorithm design reduces therefore  to the following optimization problem
with respect to $(U_t)_{t\ge 0}$, whose value we denote by $\SOC$
(stochastic optimal control)
\begin{align}
\mbox{maximize}&\;\;\;\;\;  \E\Big[\int_0^t\xi''(s) U_s\,\de s \Big]\, ,\nonumber\\
\mbox{subject to}&\;\;\;\;\; \E[M_t^2]=t \;\; \forall t\in[0,1], {\rm Law}(M_1)\in \cuP_{\Sigma}\, ,
\label{eq:SOC}\\
& \;\;\;\;\;M_t = \int_0^t \sqrt{\xi''(s)}\, U_s\, \de B_s\, .\nonumber
\end{align}
Here it is understood that $U_t$ should depend causally on the Brownian motion $(B_t)_{t\ge 0}$,
and $\cuP_{\Sigma}$ is a set of probability distributions over $\reals$
that reflect the structure of the constraint set $\Sigma_N$. Namely
\begin{align}
\Sigma_N=\S^{N-1}(\sqrt{N})\;:&\;\;\;\;\;
 \cuP_{\Sigma} = \Big\{\nu\in \cuP(\reals):\;\; \int x^2 \, \nu(\de x)=1\Big\}\, ,\\
 \Sigma_N=\{+1,-1\}^N\;:&\;\;\;\;\;
 \cuP_{\Sigma} = \Big\{\nu\in \cuP(\reals):\;\; {\rm supp}(\nu)\in[-1,+1]\Big\}\, .
 \end{align}
 Notice that naively one would have imposed the constraint ${\rm supp}(\nu)\in\{+1,-1\}$.
 However the  optimal solution to the problem with constraint $M_1\in[+1,-1]$ turns out to 
 have $M_1$ indeed concentrated on the extremes, and the interval constraint is more convenient analytically.
 
The problem \eqref{eq:SOC} is a a stochastic optimal control problem in one-dimension
and continuous time: we want to design the control process $(U_t)_{t\ge 0}$
as to satisfy the constraints and maximize the objective. It is not obvious that the sequence
of approximations (for large $N$ and small $\delta$) leading to it can indeed be implemented
as to yield an actual algorithm of the form \eqref{eq:general_amp}. However this turns out the case,
and therefore denoting by $\SOC$ the value achieved by problem \eqref{eq:SOC},
there exists an AMP algorithm such that
\begin{align}
\lim_{N\to\infty}\frac{1}{N} H_N(\bsigma^{\salg}) = \SOC\, .
\end{align}
Finally, a duality argument can be used to show that, under technical conditions 
$\SOC=\ALG$.

\subsection{Hardness}
\label{sec:Hardness}

When the infimum \eqref{eq:E-ALG} is achieved on $\gamma_{\cuL}^{\star}$ non-monotone, 
we necessarily have $\OPT>\ALG$ and therefore there is a gap between
the maximum energy and the energy achieved by the efficient algorithms described in
the previous section.
It is natural to ask whether this gap is fundamental, or there exists other polynomial-time 
algorithms that achieve a better energy than $\ALG$. 
From a physics perspective, one may wonder about the behavior of physical dynamics, such as Langevin
dynamics or simulated annealing. These are expected to converge to a threshold energy, 
which depends on the details of the dynamics (e.g. the annealing schedule)
but has proven challenging to compute (with the notable exception
of homogeneous, degree $p$ spherical models). What is the relation between $\ALG$
and the energy achieved by such annealing algorithms?

It was proven in \cite{huang2021tight} that $\ALG$ is the maximum energy achieved by a broad class of 
`Lipschitz' algorithms. This include Langevin dynamics and simulated annealing 
as special cases (run for physical time of order one as $N\to\infty$).

In order to state this result formally, we view an algorithm as a function that takes as
input an Hamiltonian and returns as output as a configuration:
\begin{align}
\cA: H_N\mapsto \cA(H_N)\in \oSigma_N\, .
\end{align}
Notice that we accept as output magnetization vectors in
$\oSigma_N=[-1,+1]^N$ (Ising) or $\oSigma_N=\Ball^N(\sqrt{N})$. 
There is no loss of generality here because magnetization vectors can be rounded to spin vectors
without decreasing significantly the energy \cite{sellke2021approximate}.

For a given Hamiltonian 
$H_N$, we denote by $\bG(H_N) = (\bG^{(k)})_{k\ge 2}$ the vector
obtained by listing all the couplings (notice that these are normalized to be typically of
order one). An algorithm $\cA$ as Lipschitz constant $L$ if, for any two Hamiltonians
$H_N$, $H_N'$:
\begin{align}
\frac{1}{\sqrt{N}}\big\| \cA(H_N) -\cA(H_N')\big\|_2\le 
L\big\|\bG(H_N)-\bG(H'_N)\big\|_2\, .
\end{align}

We can now state the Lipschitz hardness result\footnote{The statement
proven in \cite{huang2021tight} is stronger than the version given here in two directions. First,
the probability   $H_N(\bsigma^{\salg}) \ge N(\ALG+\eps)$ is shown to be exponentially
small in $N$. Second, the result holds for a somewhat broader class of ``overlap concentrated'' 
algorithms.}.
\begin{theorem}
Consider a spin glass Hamiltonian (either Ising or spherical), with 
$c_k=0$ for $k$ odd (equivalently, $\xi$ an even function). If $\cA$ is 
an $L$ Lipschitz algorithm (for $L$ a constant), then, letting $\bsigma^{\salg} =\cA(H_N)$,
we have
\begin{align}
\lim_{N\to\infty} \frac{1}{N}H_N(\bsigma^{\salg}) \le \ALG\, .
\end{align}
\end{theorem}

Remarkably, the proof of this theorem relies on a geometric
description of the computational barrier at energy $\ALG$.  Roughly speaking,
$N \cdot\ALG$ is the largest energy such that for any ultrametric
overlap structure, we can find $K$
configurations $\bsigma^1$, $\bsigma^2$, \dots $\bsigma^K$ with the assigned
energy, and whose overlap matrix matches a target overlap matrix associated to the
Parisi measure $\gamma_{\cuL}^{\star}$.

The reason why this might lead to an algorithmic barrier can be gleaned by noticing
that both algorithms described in the previous sections
can be modified to produce not one but multiple configurations.
For instance, the Hessian ascent algorithm can be randomized by using not the top 
eigenvector of the Hessian, but one of the eigenvectors with eigenvalues larger
than $(1-\eps)$ of the largest. By the semicircle law and and a large deviation bound, there
will be (with high probability) at least $c(\eps)N$ such eigenvectors for some 
constant $c(\eps)>0$. Randomized Hessian ascent follows a direction
in the subspace spanned by them uniformly at random.

 We can then generate two configurations
$\bsigma^{1,\salg}$, $\bsigma^{2,\salg}\in\S^{N-1}(\sqrt{N})$  with energies 
$H_N(\bsigma^{1,\salg}), H_N(\bsigma^{2,\salg})\approx N\ALG$ as follows.
Run two copies of the randomized Hessian ascent algorithm $\bm^{1,t}, \bm^{2,t}$
by using the same randomness for $t\le t_*$, and independent randomness (independent
choice of the Hessian eigenvector) for $t>t_*$. This results in two
spin configurations $\bsigma^{1,\salg},\bsigma^{2,\salg}$ with the claimed energy,
and $\<\bsigma^{1,\salg},\bsigma^{2,\salg}\>/N = t_*+ o_N(1)$.
Analogous constructions hold for the IAMP algorithm of Section \ref{sec:IAMP}.

This discussion clarifies why  Hessian ascent and IAMP cannot find spin configurations
at energies at which we cannot construct $K$-uples of spin configurations with arbitrary
overlap matrix. However, it does not explain why the same limitation holds
for arbitrary Lipschitz algorithms. The key idea is to run the same algorithm 
$\cA$ on two Hamiltonians $H^{(1)}_N$, $H^{(2)}_N$, that are identically distributed, jointly 
Gaussian \emph{correlated} copies of $H_N$. As a consequence 
of the Lipschitz property, the algorithms output $\cA(H^{(1)}_N)$, $\cA(H^{(2)}_N)$
have overlap that is a nearly deterministic function of the correlation 
between the Hamiltonians.

This construction is repeated in a hierarchical fashion to show that
$\cA$ can find $K$
configurations $\bsigma^1$, $\bsigma^2$, \dots $\bsigma^K$, with a certain ultrametric
overlap matrix, and a certain energy with respect to the sum of correlated Hamiltonians
$H^{(1)}_N+H^{(2)}_N+\dots+H^{(K)}_N$. Such a constellation is shown not to exist above 
energy $\ALG$ via an interpolation argument.

\subsection*{Acknowledgements}

 AA was supported by the NSF CAREER grant DMS-1653552.
AM was supported by the NSF grant CCF-2006489 and the ONR grant N00014-18-1-2729.
 Part of this work was carried out while Andrea  Montanari was on partial leave from Stanford and a 
 Chief Scientist at Ndata Inc  dba Project N. The present
research is unrelated to AM’s activity while on leave
  ES is the incumbent of the Skirball Chair in New Scientists and was supported by the Israel 
  Science Foundation
grant 2055/21.

\bibliographystyle{abbrv}

\end{document}